\documentclass[12pt]{article} \textwidth=6in \oddsidemargin=0in
\begin{document}
\def\iy{\infty} \def\be{\begin{equation}} \def\ee{\end{equation}}
\def\la{\lambda} \def\inv{^{-1}}\def\ch{\raisebox{.4ex}{$\chi$}}
\def\ov{\over} \def\s{\sigma} \def\t{\tau} \def\a{\alpha} \def\x{\xi}
\newcommand{\twotwo}[4]{\left(\begin{array}{cc}#1&#2\\&\\#3&#4\end{array}\right)}
\def\dl{\delta} \def\R{\bf R} \def\la{\lambda} \def\st{\tilde\s}
\def\spt{\widetilde{\s^+}} \def\smt{\widetilde{\s^-}}
\def\l{\lambda} \def\ov{\over} 
\def\chpm{\ch_{(-\a,\,\a)}} \def\chp{\ch_{(0,\,\a)}}
\def\chm{\ch_{(-\a,\,0)}} \def\cp{\ch^+} \def\cm{\ch^-}
\def\om{\omega} \def\d{\det\,} \def\noi{\noindent}

\hfill June 20, 2005 
\begin{center}{ \large\bf  Asymptotics of a Class of Operator Determinants}\end{center}

\begin{center}{{\bf Harold Widom}\footnote{Supported by National Science 
Foundation grant DMS-0243982.}\\
{\it Department of Mathematics\\
University of California, Santa Cruz, CA 95064\\
e-mail: widom@math.ucsc.edu}}\end{center}

\begin{abstract}
{In previous work of C. A. Tracy and the author asymptotic formulas were derived for certain operator determinants whose interest lay in the fact that quotients of them gave solutions to the cylindrical Toda equations. In the present paper we consider a more general class of operators which retain some of the properties of those cited and we find analogous asymptotics for the determinants.}
\end{abstract}
\begin{center}{\bf I. Introduction and statement of the result}\end{center}

There are innumerable instances in the mathematics and mathematical physics literature where the problem arises of determining the asymptotics of operator determinants of the form $\det\,(I+K)$ where $K$ is a trace class operator depending on a parameter. The first general result is probably the continuous analogue of the strong Szeg\"o limit theorem due to M.~Kac \cite{K}, and independently by N.~I.~Ahieser \cite{A}, where $K$ is an integral operator of convolution type on an interval whose length tends to infinity. (Thus the result is generally called the Kac-Ahieser theorem.) There have been many generalizations of this result, but by no means all such problems fall into this category and their results cannot be used directly. This paper is concerned with some of these other cases. 

In \cite{TW} the authors considered operator kernels of the form
\[K(u,v)=\int{e^{-t[(1-\om)u+(1-\om\inv)u\inv]}\over-\om u+v}\,d\rho(\om)\]
acting on $L^2(\R^+)$, where $\rho$ can be any finite complex measure supported on a compact subset of 
$\{\om\in{\bf C}:\,\Re\,\om<1,\ \Re\,\om\inv<1\}$.
The interest in these kernals was due to the fact that logarithms of  certain ratios of determinants give solutions to the cylindrical Toda equations. The main interest is in the limit $t\to0+$. In the cited paper asymptotics of the form $b\,t^a$ were determined, with the constants $a$ and $b$ having integral representations. Thus, they were explicitly determined. This was done, naturally, only under certain assumptions on $\rho$. 

If one makes the substitutions $u\to e^x,\ v\to e^y$ then operators on $L^2(\R^+)$ become operators on $L^2(\R)$. In the present paper we consider a more general class of operators on $L^2(\R)$ which retain some of the properties of those cited and we find analogous asymptotics for the determinants. The approach we use is different from that in \cite{TW} in that we deal with the determinants directly rather than through the resolvent. Because of the generality of the setting we do not obtain an explicit integral representation for the constant factor in the asymptotics; rather it itself is given in terms of operator determinants. At the end we shall indicate what special property of the particular operators of \cite{TW} allows the evaluation of these determinants.

Our setting is a family of trace class operators $K_\a$ on $L^2(\R)$ which converge strongly as $\a\to\iy$ to an operator $K$ with kernel $k(x-y)$, where $k\in L^1(\R)$. Thus each $\det\,(I+K_\a)$ is defined but $\det\,(I+K)$ is not. The problem is to find the asymptotics of $\d (I+K_\a)$. 

The first connection between $K_\a$ and $K$ we require is that
\[\ch^-\,(K_\a-K)\,\ch^+=o_1(1),\ \ \ \ch^+\,(K_\a-K)\,\ch^-=o_1(1).\]
Here $\ch^\pm$ denotes multiplication by $\ch_{\R^\pm}$ and $o_1(1)$ denotes any family of operators whose trace norms are $o(1)$. (In particular, $\ch^-\,K\,\ch^+$ and $\ch^-\,K\,\ch^+$ are trace class.)

The classical example is $K_\a=\chpm\,K\,\chpm$, whose determinant is given by the Kac-Ahieser theorem, but our operators are different. (Actually, they include these as a special case.) To state what characterizes our family $K_\a$ we introduce the translation operator $T_a$ defined by $T_a f(x)=f(x-a)$. Our main assumption is that there are operators $K_\pm$ such that 
\[K_+-\ch^-\,K\,\ch^-\ \ {\rm and}\ \ \ K_--\ch^+\,K\,\ch^+\]
are trace class and 
\be\ch^+(\,K_a-T_\a\,K_+\,T_{-\a})\,\ch^+=o_1(1),\ \ \ 
\ch^-\,(K_a-T_{-\a}\,K_-\,T_{\a})\,\ch^-=o_1(1).\label{cond1}\ee

A basic requirement is that if $\hat k$ is the Fourier transform of $k$,
\[\hat k(\x)=\int_{-\iy}^\iy e^{ix\x}\,k(x)\,dx,\]
and if we set $\s(\x)=1+\hat k(\x)$ then
\be \s(\x)\ne0,\ \ \ {\rm arg}\;\s(\x)|_{-\iy}^\iy=0.\label{index}\ee
This assures that the Wiener-Hopf operators $I+\ch^\pm\,K\,\ch^\pm$ are invertible.\footnote{See, for example, \cite[\S I.8]{GF}.} We shall also assume that $|x|^{1/2}\,k(x)\in L^2(\R)$ which, together with (\ref{index}), will allow us to use the Kac-Ahieser theorem.\footnote{This says that if $W_\a(\s)=I+\chpm\,K\,\chpm$ acting on $L^2(-\a,\a)$ then 
$\d W_\a(\s)\sim G(\s)^{2\a}\,E(\s)$, where $G(\s)=\exp\,(s(0)),\ E(\s)=\exp\,\{\int_0^\iy x\,s(x)\,s(-x)\,dx\}$, and $s(x)$ is the inverse Fourier transform of $\log \s(\x)$. There is also the alternative expression $E(\s)=\d (W(\s)\,W(\s\inv))$.}

Our result is that
\be\det(I+K_\a)\sim G(\s)^{2\a}\,E(\s)\,\d(I+K_1)\,\d(I+K_2),\label{result}\ee
where $G(\s)$ and $E(\s)$ are the constants in the Kac-Ahieser theorem and $K_1$ and $K_2$ are certain trace class operator on $L^2(\R)$ built out of $K$ and $K_\pm$. Precisely,
\[K_1=(I+\ch^+\,K\ch^+)\inv\,(K_--\ch^+\,K\,\ch^+),\ \ \ K_2=
(I+\ch^-\,K\ch^-)\inv\,(K_+-\ch^-\,K\,\ch^-).\]

In the classical example where $K_\a=\chpm\,K\,\chpm$ we take $K_+=\ch^-\,K\,\ch^-$ and $K_-=\ch^+\,K\,\ch^+$  (the operators in (\ref{cond1}) are then identically zero) so $K_1=K_2=0$. 

In the example which arises from the simplest case in \cite{TW} (after symmetrization and variable change) $K_\a$ is the integral operator with kernel
\[K_\a(x,y)=\la \,{e^{-(e^{x-\a}+e^{y-\a}+
e^{-x-\a}+e^{-y-\a})}\over{\rm cosh}\,[(x-y)/2]},\]
and the other kernels are
\[K(x,y)=\la \;{1\over{\rm cosh}\,[(x-y)/2]},\]
\[K_+(x,y)=\la \;{e^{-(e^{x}+e^{y})}\over{\rm cosh}\,[(x-y)/2]},\]
\[K_-(x,y)=\la \;{e^{-(e^{-x}+e^{-y})}\over{\rm cosh}\,[(x-y)/2]}.\]
The determinants on the right side of (\ref{result}) can be evaluated in this case.

\noi{\bf Remark}. The cases where the conditions (\ref{index}) are satisfied are the easier ones. There is great interest in the integrable systems community in operators for which $\s(\x)$ may have zeros. (In the above example this occurs when $\l\in(-\iy,-\pi\inv]$.) Asymptotic formulas have been derived in some of these cases \cite{MTW,Ki}, but as far as we know none have been proved rigorously. A hope is that the approach presented here may be applicable, at least to some extent, to these. 
\pagebreak

\begin{center}{\bf II. Derivation of the result}\end{center}

We think of $L^2(\R)$ as $L^2(\R^-)\oplus L^2(\R^+)$ and
the corresponding matrix representation. Condition (\ref{cond1}) tells us that with error $o_1(1)$ the matrix representation of 
$K_\a$ is 
\[\twotwo{\cm\,T_{-\a}\,K_-\,T_{\a}\,\cm}{K}{K}{\cp\,T_\a\,K_+\,T_{-\a}\,\cp}.\]

We define
\[K_{11}=K_--\ch^+\,K\,\ch^+\ \ {\rm and}\ \ \ K_{22}=K_+-\ch^-\,K\,\ch^-,\]
which by assumption are trace class operators on $L^2(\R)$. The upper-left corner of the above matrix may be written
\be \cm\,(T_{-\a}\,\ch^+\,K\,\ch^+\,T_{\a}+T_{-\a}\,K_{11}\,T_{\a})\,\cm.\label{ul}\ee
We shall use, here and below, the fact that $K$ commutes with translations and the relation
\be\ch_J\,T_a=T_a\,\ch_{J-a}\label{JT}\ee
for any $a$ and any set $J$. Thus (\ref{ul}) may be written $\chm\,K\,\chm+\cm\,T_{-\a}\,K_{11}\,T_{\a}\,\cm$. Similarly for the lower-right corner, and so the matrix representation of $K_\a$ is $o_1(1)$ plus
\[\twotwo{\chm\,K\,\chm+\cm\,T_{-\a}\,K_{11}\,T_{\a}\,\cm}{\cm\,K\,\cp}{\cp\,K\,\cm}
{\hspace{-2em}\chp\,K\,\chp+\cp\,T_{\a}\,K_{22}\,T_{-\a}\,\cp}.\]

The upper-right corner is trace class, so if we multiply it on either side by $\ch_{(-\iy,\,-\a)}$ or $\ch_{(\a,\,\iy)}$ the result is $o_1(1)$. Similarly for the lower-left corner. Hence with this error the above equals
\[\twotwo{\chm\,K\,\chm}{\chm\,K\,\chp}{\chp\,K\,\chm}
{\chp\,K\,\chp}+
\twotwo{\cm\,T_{-\a}\,K_{11}\,T_{\a}\,\cm}{0}{0}
{\cp\,T_{\a}\,K_{22}\,T_{-\a}\,\cp}.\]

The operator
\[I+\twotwo{\chm\,K\,\chm}{\chm\,K\,\chp}{\chp\,K\,\chm}
{\chp\,K\,\chp}\]
is just $W_\a(\s)$ in its matrix representation and so its determinant is asymptotically equal to $G(\s)^{2\a}\,E(\s)$, by the Kac-Ahieser theorem.\footnote{Strictly speaking $W_\a(\s)$ acts on $L^2(-\a,\a)$ while the above operator acts on $L^2(\R)$. But the determinants are the same.}

The next step is to factor out this operator from $I+K_\a$
and determine the asymptotics of the determinant of the result. For convenience we write the above operator as $D+E$, where
\[D=\twotwo{I+\chm\,K\,\chm}{0}{0}{I+\chp\,K\,\chp}\]
and
\[E=\twotwo{0}{\chm\,K\,\chp}{\chp\,K\,\chm}{0},\]
so the operator whose determinant we now want is $o_1(1)$ plus\footnote
{The reason the error remains $o_1(1)$ after multiplying by $W_\a(\s)\inv$ is that these operators have uniformly bounded norms.
See \cite[\S III.1]{GF}.}
\be I+(D+E)\inv\,\twotwo{\cm\,T_{-\a}\,K_{11}\,T_{\a}\,\cm}{0}{0}
{\cp\,T_{\a}\,K_{22}\,T_{-\a}\,\cp}.\label{IDE}\ee
In our notation $D$ stands for ``diagonal'' and $E$ stands for ``error'' because, as we shall now show, its contribution here will be $o(1)$.

The difference between the above operator and the one with $E$ replaced by 0 equals
\[(D+E)\inv\,E\,D\inv\,\twotwo{\cm\,T_{-\a}\,K_{11}\,T_{\a}\,\cm}{0}{0}
{\cp\,T_{\a}\,K_{22}\,T_{-\a}\,\cp}.\]
The left-most operator $(D+E)\inv$ is $W_\a(\s)\inv$, which has uniformy bounded norm. We shall show that the product of the remaining ones is $o_1(1)$.

We use the notations
\[W_\a^+=I+\chp\,K\,\chp,\ \ \ W_\a^-=I+\chm\,K\,\chm\]
and
\[W^+=I+\ch^+\,K\,\ch^+,\ \ \ W^-=I+\ch^-\,K\ch^-\]
because these operators arise so often. Observe that by (\ref{JT})
\be W_\a^+\,T_\a=T_\a\,W_\a^-,\ \ \ (W_\a^+)\inv\,T_\a=T_\a\,(W_\a^-)\inv.\label{WT}\ee

The upper-right corner of the product in question is \[\chm\,K\,\chp\,(W_\a^+)\inv\,\ch^+\,T_{\a}\,K_{22}\,T_{-\a}\,\ch^+.\]
The last $T_{-\a}\,\ch^+$ may be ignored since its norm is $O(1)$. By (\ref{JT}) and (\ref{WT}) what remains is equal to
\be\chm\,K\,\chp\,T_\a\,(W_\a^-)\inv\,\ch_{(-\a,\iy)}\,K_{22}.\label{entry}\ee
The operators $(W_\a^-)\inv$ converge strongly to $(W^-)\inv$.\footnote{For this also see \cite[\S III.1]{GF}.} Since $K_{22}$ is trace class it follows that $(W_\a^-)\inv\,\,\ch_{(-\a,\iy)}\,K_{22}$ converges in trace norm. As for the remaining factors, we observe that $\chm\,K\,\chp$ converges in norm to the compact operator $\ch^-\,K\,\ch^+$ while $T_\a$ converges weakly to zero. It follows that $\chm\,K\,\chp\,T_\a$ converges strongly to zero, and we deduce that (\ref{entry}) converges to zero in trace norm. A similar argument applies to the lower-left corner of the product.

We have shown that with error $o_1(1)$ we may replace $E$ by zero in (\ref{IDE}), which then becomes 
\[I+\twotwo{\ch^-\,(W_\a^-)\inv\,
T_{-\a}\,K_{11}\,T_{\a}\,\ch^-}{0}{0}
{\hspace{-2em}\ch^+\,(W_\a^+)\inv\,T_{\a}\,K_{22}\,T_{-\a}\,\ch^+}.\]
This operator, which  acts on $L^2(\R^-)\oplus L^2(\R^+)$, can be extended in an obvious way, without change of notation or determinant, to one acting on $L^2(\R)\oplus L^2(\R)$. By (\ref{WT}) and (\ref{JT}) we may now rewrite it as
\[I+\twotwo{T_{-\a}\,\ch_{(-\iy,\a)}\,(W_\a^+)\inv\,
K_{11}\,\ch_{(-\iy,\a)}\,T_{\a}}{0}{0}
{\hspace{-4em}T_{\a}\,\ch_{(-\a,\iy)}\,
(W_\a^-)\inv\,K_{22}\,\ch_{(-\a,\iy)}\,T_{-\a}}.\]
If we left-multiply by the unitary operator $\left(\begin{array}{cc}
T_{\a}&0\\0&T_{-\a}\end{array}\right)$ and right-multiply by its inverse $\left(\begin{array}{cc}
T_{-\a}&0\\0&T_{\a}\end{array}\right)$ the operator becomes
\[I+\twotwo{\ch_{(-\iy,\a)}\,(W_\a^+)\inv\,
K_{11}\,\,\ch_{(-\iy,\a)}}{0}{0}
{\hspace{-2em}\ch_{(-\a,\iy)}\,(W_\a^-)\inv\,K_{22}\,\ch_{(-\iy,\a)}}.\]
The determinant is unchanged, and the error term remains $o_1(1)$. 
By an argument already used, this operator (plus the error term $o_1(1)$) converges in trace norm to 
\[I+\twotwo{(W^+)\inv\,
K_{11}}{0}{0}
{\hspace{-1em}(W^-)\inv\,K_{22}},\]
and so its determinant converges to the determinant of this one. The determinant of this one equals the product of determinants on the right side of (\ref{result}), which is therefore now established.

\begin{center}{\bf III. Final remarks}\end{center}

\noi{\bf Remark 1}. There is the possibility that $I+K_1$ or $I+K_2$ is not invertible, in which case the product of determinants on the right side of (\ref{result}) is zero. Since
\[I+K_1=(I+\ch^+\,K\,\ch^+)\inv\,(I+K_-),\]
the invertibility of $I+K_1$ is equivalent to that of $I+K_-$. Similarly the invertibility of $I+K_2$ is equivalent to that of $I+K_+$. These are separate issues.

\noi{\bf Remark 2}. Here is a rough explanation of why the determinants on the right side of (\ref{result}) are sometimes evaluable. The first one, for example, is
\[\d[I+(I+\ch^+\,K\,\ch^+)\inv\,K_{11}].\]
If we introduce a paramenter $\l$ then
\[{d\over d\l}\log\,\d[I+(I+\l\,\ch^+\,K\,\ch^+)\inv \l\,K_{11}]\]
is equal to the trace of
\[(I+\l\,K_-)\inv\,K_--(I+\l\,\ch^+\,K\,\ch^+)\inv\,\ch^+\,K\,\ch^+.\]
Now $(I+\l\,\ch^+\,K\,\ch^+)\inv$ is the direct sum of $I$ acting on 
$L^2(\R^-)$ and $W(1+\l\,\hat k)\inv$ acting on $L^2(\R^+)$, so it is known --- it is expressible in terms of the Wiener-Hopf factors of $1+\l\,\hat k$. 
The pleasant fact is that the first operator, while not of this form, can be brought to this form in the cases considered in \cite{TW}. 

In the particular one mentioned in the introduction the kernel of $K_-$ equals (when $\l$ there is taken to be 1)
\[2\,e^{-(x+y)/2}\,{e^{-(e^{-x}+e^{-y})}\ov e^x+e^y}.\]
This has the integral representation
\[2\,e^{-(x+y)/2}\int_0^\iy e^{-({e^{u-x}+e^{u-y}})}\,e^u\,du.\]
Hence if $M$ is the integral operator from $L^2(\R^+)$ to $L^2(\R)$ with kernel
\[M(x,u)=\sqrt{2}\,e^{-x/2} \,e^{-{e^{u-x}}}\,e^{u/2}\]
and $N$ is the integral operator from $L^2(\R)$ to $L^2(\R^+)$ with kernel 
\[N(v,y)=\sqrt{2}\,e^{-y/2} \,e^{-{e^{v-y}}}\,e^{v/2}\]
then $K_-=MN$. The kernel of $NM$, an operator on $L^2(\R^+)$, is
\[\int _{-\iy}^\iy N(u,x)\,M(x,v)\,dx=
2\,{e^{(u+v)/2}\ov e^u+e^v}={1\ov \cosh[(u-v)/2]}.\]
This is itself a Wiener-Hopf operator so the inverse of $I+\l\,NM$ may be written down, and we have
\[(I+\l\,K_-)\inv\,\l\, K_-=I-(I+\l\,MN)\inv=\l\,M\,(I+\l\,NM)\inv\,N.\]
This enables one, at least in principle, to compute the logarithmic derivative of the determinant and hence the determinant itself.

\end{document}